\documentclass{article}
\usepackage{amsmath, amssymb}
\usepackage{pgfplots}
\pgfplotsset{compat=1.18}
\usepackage{theorem}

\usepackage[hmargin=2.5cm, vmargin=1.1in]{geometry}
\usepackage{parskip}

\numberwithin{equation}{section}
{\theoremstyle{break}\newtheorem{theorem}{Theorem}[subsection]}
{\theoremstyle{break}\theorembodyfont{\rmfamily}\newtheorem{example}[theorem]{Example}}
{\theoremstyle{break}\theorembodyfont{\rmfamily}\newtheorem{definition}[theorem]{Definition}}
{\theoremstyle{break}\theorembodyfont{\rmfamily}\newtheorem{notation}[theorem]{Notation}}
{\theoremstyle{break}\theorembodyfont{\rmfamily}\newtheorem{remark}[theorem]{Remark}}

\newenvironment{proof}[1][Proof]{\noindent\textbf{#1.} }{\
\rule{0.5em}{0.5em}}

\allowdisplaybreaks[2]

\title{Necessary And Sufficient Characterization Of Absolutely Continuous Functions Defined Over Unbounded Intervals.}
\author{Gourav Banerjee}

\begin{document}

\maketitle

\begin{abstract}
In this paper, we investigate and find a necessary and sufficient condition for a function to be absolutely continuous over $\mathbb{R}$ (denoted by $AC(\mathbb{R})$) or any unbounded interval in $\mathbb{R}$ . Note that the Lebesgue's Fundamental theorem of Calculus gives us a necessary and sufficient condition\cite{book:B} for a function defined over a closed interval [a,b] to be absolutely continuous ,and the condition is that the derivative of the function should be in $L^1_{loc}([a,b])$. However, we don't have any such sufficient condition on the derivative of a function that is absolutely continuous over unbounded intervals.
One necessary condition is that the function must be locally absolutely continuous (denoted by $AC_{loc}(\mathbb{R})$), but it may not be globally absolutely continuous despite being locally absolutely continuous(we give an explicit example of this).
\\
The theorem 1 in this paper gives us a necessary and sufficient condition for a function belonging in $AC_{loc}(\mathbb{R})$ to belong in $AC(\mathbb{R})$ in terms of its derivative and identifies the space to which the derivative of an $AC(\mathbb{R})$ function must belong to as $L^1_{G}(\mathbb{R})$ (a strict subspace of $L^1(\mathbb{R})$).
\\
Moreover, we define a new space of functions called $L^1_{H}(\mathbb{R})$, and in theorem 2 we show that $L^1_G \subset L^1_H$, which helps us to find an easier criteria to check whether a function belonging to $AC_{loc}(\mathbb{R})$, belongs to $AC(\mathbb{R})$ or not.
\\
Finally, we provide a Venn diagram to explicitly show the relation of the newly defined spaces $L^1_{G}(\mathbb{R})$ and $L^1_{H}(\mathbb{R})$ with respect to the spaces $L^1_{loc}(\mathbb{R})$, $L^1(\mathbb{R})$ and $L^\infty(\mathbb{R})$.    
\end{abstract}

\tableofcontents

\section{Introduction}

The classical Fundamental Theorem of Calculus establishes a deep relationship between
differentiation and integration, and it lies at the foundation of real analysis. In the
Riemann setting, the theorem guarantees that if $f$ is continuous on a compact interval,
then the function $F(x) = \int_a^x f(t)\, dt$ is differentiable and $F' = f$. However,
Riemann integration is too restrictive to capture many naturally occurring functions,
leading to the development of the Lebesgue integral and the modern formulation of the
Fundamental Theorem of Calculus.

Lebesgue's Second Fundamental Theorem of Calculus provides a complete characterization of
absolutely continuous functions on compact intervals. Specifically, a function $F$ defined
on $[a,b]$ is absolutely continuous if and only if $F'$ exists almost everywhere, belongs
to $L^1([a,b])$, and
\[
F(x) = F(a) + \int_a^x F'(t)\, dt.
\]
Thus, on bounded intervals, absolute continuity is precisely the condition ensuring that
differentiation and integration invert each other in the Lebesgue sense. This identifies
$AC([a,b])$ as the exact class of functions representable as indefinite integrals of
integrable functions.

However, when the domain is unbounded (e.g., $\mathbb{R}$ or $(a,\infty)$), this
characterization no longer applies directly. A function may be absolutely continuous on
every compact subinterval yet fail to be absolutely continuous globally, so that
\[
AC(\mathbb{R}) \subsetneq AC_{loc}(\mathbb{R}).
\]
For instance, $f(x)=\sqrt{|x|}$ extended periodically is locally absolutely continuous but
not absolutely continuous on $\mathbb{R}$. Moreover, the condition $F' \in L^1(\mathbb{R})$
is too strong, since even the simple function $f(x)=x$ belongs to $AC(\mathbb{R})$ but
satisfies $f' \notin L^1(\mathbb{R})$. This shows that global absolute continuity depends on
a subtler form of integrability for the derivative.

To capture this, we introduce the function space $L^1_G(\mathbb{R})$, consisting of
measurable functions whose integrals are finite over all measurable subsets of finite
measure. Our first main result shows that this space provides the exact integrability
criterion needed for global absolute continuity:

\[
f \in AC(\mathbb{R}) \quad \Longleftrightarrow \quad f \in AC_{loc}(\mathbb{R}) \text{ and } f' \in L^1_G(\mathbb{R}).
\]

This theorem extends Lebesgue’s characterization from compact intervals to unbounded
domains, thereby establishing $L^1_G(\mathbb{R})$ as the natural analogue of $L^1([a,b])$
for the global setting.

We also introduce a second function space, $L^1_H(\mathbb{R})$, defined in terms of the
measure of the level sets of $|f|$. In Theorem 2, we show that
\[
L^1_G(\mathbb{R}) \subset L^1_H(\mathbb{R}),
\]
and that the inclusion is proper. This provides an alternative, often more accessible
condition for verifying global absolute continuity. The relationships between these spaces
and the classical Lebesgue spaces are illustrated in a diagram later in the paper.

Together, Theorem 1 and Theorem 2 establish a structural framework for understanding
absolute continuity over unbounded intervals and clarify the precise role of derivative
integrability in the global setting.

\section{Function spaces of our interest}
We begin this section by setting up the notation for the function spaces that will be mostly used throughout this paper.

\begin{definition}[Absolutely Continuous Function]
A function $f : \mathbf{I} \to \mathbf{R}$ ($\mathbf{I}$ being an interval in $\mathbf{R}$) is said to be \emph{absolutely continuous} on $\mathbf{I}$ if for every $\varepsilon > 0$, there exists a $\delta > 0$ such that for any finite collection of pairwise disjoint subintervals $\{(x_k, y_k)\}_{k=1}^n$ of $\mathbf{I}$ satisfying
\[
\sum_{k=1}^n (y_k - x_k) < \delta,
\]
it follows that
\[
\sum_{k=1}^n |f(y_k) - f(x_k)| < \varepsilon.
\]
\end{definition}

\begin{notation}
We define\\
$AC(\mathbf{I})$= $\{$space of measurable functions  f: $\mathbf{I} \to \mathbf{R}$ that are absolutely continuous over $\mathbf{I}$ $\}$.
\end{notation}

\begin{notation}[Locally Absolutely Continuous Function]
We define\\
$AC_{loc}(\mathbf{I})$ = $\{$ space of measurable functions f: $\mathbf{I} \to \mathbf{R}$ that are absolutely continuous over any compact sub-interval of $\mathbf{I}$ .$\}$
\end{notation}

Now we introduce two new spaces of functions which are an integral part of the main result of this paper.

\begin{notation}
We define\\
$\mathbf{L^1_G(\mathbf{I})}$ = $\{$space of measurable functions f: $\mathbf{I} \to \mathbf{R}$ such that $ \int_K |f| < \infty$
 where K is a Lebesgue measurable subset of I and has finite measure.$\}$
\end{notation}

\begin{notation}
We define\\
$\mathbf{L^1_H(\mathbf{I})}$ =$\{$set of all measurable functions $f : \mathbf{I} \to \mathbf{R}$ such that  for a given $f$ there exists $ M > 0 : \mu\{x \in \mathbf{I} \mid |f(x)| \geq M\} <\infty \}$.
\end{notation}
Now we provide a Venn diagram to explicitly show the relationship between these newly defined spaces and the spaces {$L^1(\mathbb{R})$}, {$L_{loc}^1(\mathbb{R})$, {$L^\infty(\mathbb{R})$  with some examples.

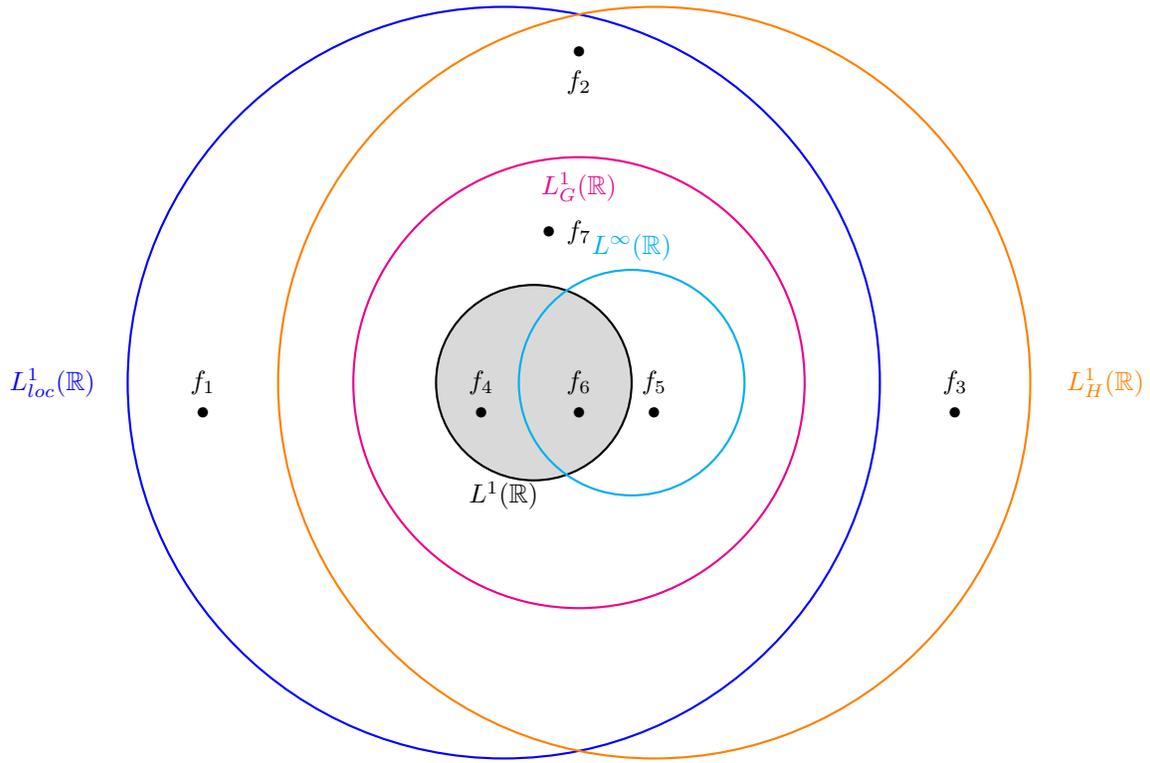
\begin{center}
\begin{tikzpicture}
    \fill[gray!30] (-0.6,0) ellipse (1.3 and 1.3);
    \draw[thick] (-0.6,0) ellipse (1.3 and 1.3);
    \draw[thick,color=blue] (-1,0) ellipse (5 and 5);
    \draw[thick,color=orange] (1,0) ellipse (5 and 5);
    \draw[thick,color=magenta] (0,0) ellipse (3 and 3);
    \draw[thick,color=cyan] (0.7,0) ellipse (1.5 and 1.5);

    \node at (-1, -1.5) {$L^1(\mathbb{R})$};
    \node[color=magenta] at (0, 2.6) {$L_{G}^1(\mathbb{R})$};
    \node[color=blue] at (-7, 0) {$L_{loc}^1(\mathbb{R})$};
    \node[color=cyan] at (0.7, 1.8) {$L^\infty(\mathbb{R})$};
    \node[color=orange] at (7, 0) {$L_{H}^1(\mathbb{R})$};
    \node at (-5, 0) {$f_1$};
    \node at (0, 4) {$f_2$};
    \node at (5, 0) {$f_3$};
    \node at (-1.3, 0) {$f_4$};
    \node at (1, 0) {$f_5$};
    \node at (0, 0) {$f_6$};
    \node at (0, 2) {$f_7$};
    \node at (1, -0.4) {\textbullet};
    \node at (0, -0.4) {\textbullet};
    \node at (-1.3, -0.4) {\textbullet};
    \node at (-0.4, 2) {\textbullet};
    \node at (5, -0.4) {\textbullet};
    \node at (-5, -0.4) {\textbullet};
    \node at (0,4.4) {\textbullet};
    
\end{tikzpicture}
\end{center}

\bigskip

\begin{align*}
    f_1(x) &= x \\
    f_2(x) &= \sum_{k=1}^{\infty} n \cdot \chi_{[n, n + \frac{1}{n^2}]} \\
    f_3(x) &= 
    \begin{cases}
     \dfrac{1}{x}, x\ne0\\
     0,  \text{otherwise}
    \end{cases}\\
    f_4(x) &=
    \begin{cases}
    \dfrac{1}{\sqrt{|x|}}, & -1 \le x \le 1, \text{and } x\ne0 \\[6pt]
    0, & \text{otherwise.}
    \end{cases} \\
    f_5(x) &= 1 \quad \text{(constant function)} \\
    f_6(x) &=
    \begin{cases}
    1, & -1 \le x \le 1, \\[6pt]
    0, & \text{otherwise.}
    \end{cases} \\
    f_7(x) &=
    \begin{cases}
    \dfrac{1}{\sqrt{|x|}}, & -1 \le x \le 1, \text{and } x\ne0 \\[6pt]
    1, & \text{otherwise.}
    \end{cases}
\end{align*}

\begin{remark}
Here we give a brief justification for the location of the function $f_1(x)=x$ in our Venn diagram.

Consider the set 
\[
S = \bigcup_{n=1}^\infty [n, n + \tfrac{1}{n^2})
\]
Note that 
\[
m(S) = \sum_{n=1}^\infty \frac{1}{n^2} < \infty
\]
but 
\[
\int_S f_1 \, dm = \sum_{n=1}^\infty n \cdot \frac{1}{n^2} = \sum_{n=1}^\infty \frac{1}{n} = \infty
\]

\[
\therefore \quad f_1 \notin L^1(\mathbb{R}), \text{ but } f_1 \in L^1_{\text{loc}}(\mathbb{R}) 
\]
as \( f_1 \) is continuous and bounded over any compact set.

Also, given the definition of \( L^1_H(\mathbb{R}) \), it’s clear that \( f_1(x) \notin L^1_H(\mathbb{R}) \).
\end{remark}

The only thing that remains to be proved here is that $L^1_G(\mathbb{R}) $is a proper subset of$L^1_H(\mathbb{R})$ which we prove in the next section in Theorem 3.0.5.

\section{Fundamental Theorem Of Calculus}

Over the development of real analysis, various notions of integration have been introduced to extend the applicability of the Fundamental Theorem of Calculus (FTC) beyond the limitations of the classical Riemann integral. Each successive definition generalizes the previous one, allowing a broader class of functions to be integrated and differentiated while preserving, in some form, the connection between differentiation and integration. Below, we briefly outline the evolution of these integrals.

\paragraph{1. Riemann Integral.}
The Riemann integral, introduced in the 19th century, defines the integral of a function as the limit of Riemann sums over partitions of an interval. It works well for continuous or piecewise continuous functions defined on compact intervals. The classical \textit{Fundamental Theorem of Calculus} for the Riemann integral states that if \( f \) is continuous on \([a,b]\) and \( F(x) = \int_a^x f(t)\,dt \), then \( F'(x) = f(x) \). Conversely, if \( F \) is differentiable with continuous derivative \( F' \), then \( \int_a^b F'(x)\,dx = F(b) - F(a) \).

\paragraph{2. Lebesgue Integral.}
Henri Lebesgue extended the notion of integration by measuring the ``size'' of sets where a function takes certain values rather than subdividing the domain into intervals. This allowed integration of functions with far more complex sets of discontinuities. The \textit{Lebesgue's Fundamental Theorem of Calculus} (sometimes called the \textit{Second Fundamental Theorem of Calculus}) states that a function \( F:[a,b]\to\mathbb{R} \) is absolutely continuous if and only if \( F' \in L^1([a,b]) \) and
\[
F(x) = F(a) + \int_a^x F'(t)\,dt.
\]
Thus, absolute continuity exactly characterizes those functions for which the FTC holds in the Lebesgue sense.

\paragraph{3. Henstock–Kurzweil (HK) Integral.}
The Henstock–Kurzweil integral (also called the \textit{generalized Riemann} or \textit{gauge integral}) further extends integration to functions that may fail to be Lebesgue integrable, while still satisfying the FTC. Every derivative is HK-integrable, and for any function \( f \), if \( f' \) exists almost everywhere, then
\[
\int_a^b f'(x)\,dx = f(b) - f(a)
\]
in the Henstock–Kurzweil sense. The HK integral thus restores the full converse of the FTC.


\bigskip

\noindent
In this paper, however, we focus on the \textbf{Lebesgue’s Second Fundamental Theorem of Calculus}, as it provides a necessary and sufficient condition for a function defined over a \textit{compact interval} \([a,b]\) to be absolutely continuous—namely, that its derivative belongs to \(L^1([a,b])\). Nevertheless, this characterization fails to extend directly to functions defined on \textit{unbounded intervals}. The primary objective of this work is to bridge this gap by developing an analogous necessary and sufficient condition for absolute continuity on unbounded domains in \(\mathbb{R}\).

Let $f : [a,b] \to \mathbb{R}$ be an absolutely continuous function. Then:
\begin{enumerate}
    \item $f$ is differentiable almost everywhere on $[a,b]$,
    \item its derivative $f'$ is Lebesgue integrable on $[a,b]$, and
    \item for all $x \in [a,b]$,
    \[
    f(x) = f(a) + \int_a^x f'(t)\, dt.
    \]
\end{enumerate}
Conversely, if $g \in L^1([a,b])$ and we define
\[
F(x) = \int_a^x g(t)\, dt,
\]
then $F$ is absolutely continuous on $[a,b]$, $F' = g$ almost everywhere, and
\[
F(x) = F(a) + \int_a^x F'(t)\, dt.
\]
\\The theorem shows that \emph{a function is absolutely continuous on $[a,b]$ (bounded interval) if and only if it is the indefinite integral of some Lebesgue integrable function.} Thus, absolute continuity provides the exact class of functions for which the classical Fundamental Theorem of Calculus extends to the Lebesgue integral setting.

\vspace{0.25cm}However, we don't have any such sufficient condition on the derivative of a function which is in $AC(\mathbb{R})$. For example,
$f(x)=x \in AC(\mathbb{R})$ but $f'(x) \notin L^1(\mathbb{R})$. One necessary condition on the derivative of a function $f$ belonging to $AC(\mathbb{R})$ can be easily noted by noting that such a function must belong to$AC_{loc}(\mathbb{R})$, and hence $f' \in L^1_{loc}(\mathbb{R})$ but the converse is not true. We give an example\cite{book:A} of such a function which is in $AC_{loc}(\mathbb{R})$ and not in $AC(\mathbb{R})$ .\\

\begin{example}\label{ex1}
Let $f(x)=\sqrt{|x|}, -1 \leqslant x \leqslant 1$.  Let $f(x+2k)=f(x)$ for each $k \in Z$ and each $x$. 

From its definition, $f$ is continuous on $[0,1]$ and so on $\mathbb{R}$.\\
Given $\delta$, $0<\delta<1 / 2$, \\let $x_i=2i, y_i=2i+\delta / i^2$ \\ Then, for each $n$, 
\[
 \sum_{i=1}^n\left|x_i-y_i\right|<2 \delta
\] 
but 
\[
\sum_{i=1}^n\left|f\left(x_i\right)-f\left(y_i\right)\right|=\sum_{i=1}^n \frac{\sqrt{ \delta}}{i}
\] 
which tends to infinity with $n$ tending to infinity. So $f$ is not absolutely continuous.
\end{example}

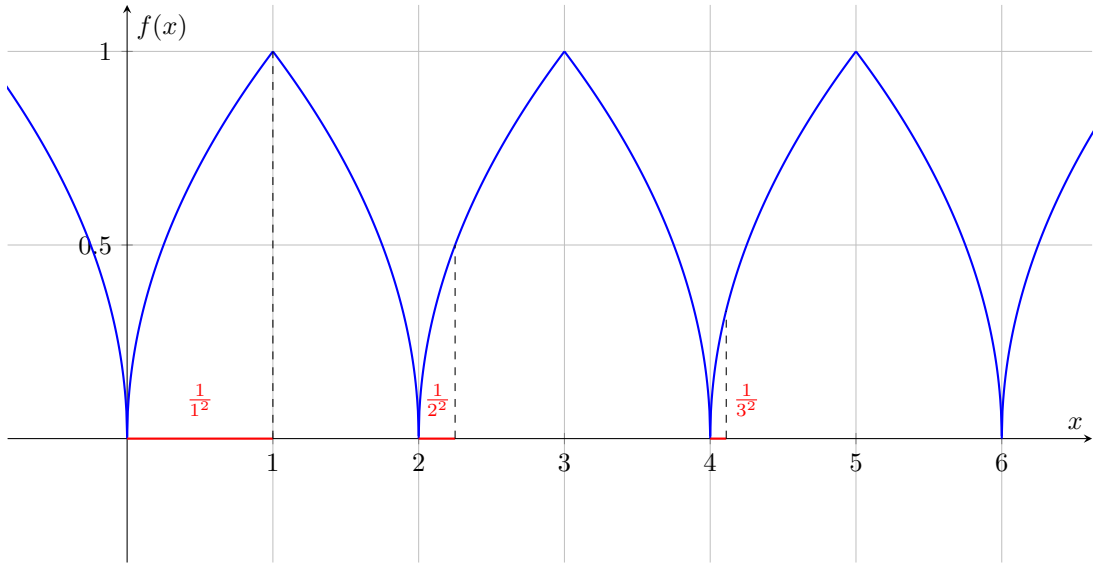
\begin{tikzpicture}
    \begin{axis}[
        axis lines=middle,
        xlabel={$x$},
        ylabel={$f(x)$},
        xtick={-2, -1, 0, 1, 2, 3, 4, 5, 6,7,8,9,10,11,12},
        ytick={0, 0.5, 1},
        ymin=-0.2, ymax=1,
        xmin=-0.2, xmax=6,
        domain=-1:6,
        samples=200,
        width=16cm,
        height=9cm,
        grid=both,
        major grid style={line width=.2pt,draw=gray!50},
        minor grid style={line width=.1pt,draw=gray!20},
        enlargelimits
    ]
        
        
        \foreach \k in {-2,0,2,4,6,8,10,12} {
            \addplot[blue, thick, domain=\k:\k+1] {sqrt(x - \k)};
        }     
        \foreach \k in {0,2,4,6,8,10} {
            \addplot[blue, thick, domain=\k-1:\k] {sqrt(-x + \k)};
        }
        \addplot[red, thick, domain=0:1] {0};
        \addplot[red, thick, domain=2:2.25] {0};
        \addplot[red, thick, domain=4:4.11] {0};
         \draw [dashed] (1,0) -- (1,1);
         \draw [dashed] (2.25,0) -- (2.25,0.5);
         \draw [dashed] (4.11,0) -- (4.11,0.33);
  \node [color=red] at (0.5,0.1) {$\frac{1}{1^2}$};
  \node [color=red] at (2.13,0.1) {$\frac{1}{2^2}$};
  \node [color=red] at (4.25,0.1) {$\frac{1}{3^2}$};

    \end{axis}
\end{tikzpicture}

Hence, the derivative of a function $f \in AC_(\mathbb{R})$ must belong to a smaller subset of $L^1_{loc}(\mathbb{R})$. It turns out that this subset is rather a strict subspace of $L^1_{loc}(\mathbb{R})$ and we call it $L^1_{G}(\mathbb{R})$.
Thus, theorem 1 gives us a necessary and sufficient condition for a function belonging in $AC_{loc}(\mathbb{R})$ to belong in $AC(\mathbb{R})$ in terms of its derivative.
Moreover, we define a new space of functions called $L^1_{H}(\mathbb{R})$, and in theorem 2 we show that $L^1_G \subset L^1_H$, which helps us to find an easier criteria to check whether a function belonging to $AC_{loc}(\mathbb{R})$, belongs to $AC(\mathbb{R})$ or not.\\

Now we present the main result of this paper.

\begin{theorem}
A locally absolutely continuous function $f$ is globally absolutely continuous if and only if $f'$ lies in $L^1_G(\mathbb{R})$  \\i.e., $f\in AC_{loc}(\mathbf{R})$ and $f' \in L^1_G(\mathbf{R}) \iff f\in AC(\mathbb{R}).$ \\
Here $L^1_G(\mathbb{R})$ represents the space of measurable functions over $\mathbb{R}$ such that their integral over every finite measurable set is finite.
\end{theorem}

\begin{remark}
A locally absolutely continuous function $f$ is almost everywhere differentiable in $\mathbf{R}$, hence we can talk about its derivative $f'$.
\end{remark}

 \begin{proof}
Let us assume that $f\notin AC(\mathbb{R})$,but$f\in AC_{loc}(\mathbb{R}) $ and  $f' \in L^1_G(\mathbb{R})$, 

Then, $\exists \, \epsilon > 0$, such that for any $\delta > 0$, there exists a sequence of disjoint intervals $(a_i, b_i), \, i \in I$ such that:
\begin{equation}\label{r1}
\sum_{i=1}^{\infty} (b_i - a_i) < \delta \quad \text{but} \quad \sum_{i=1}^{\infty} \left| f(b_i)-f(a_i)  \right| \geq \epsilon.
\end{equation}

Consider a positive sequence $(r_n)$ such that:
\[
\sum_{n=1}^{\infty} r_n < \infty. \quad \text{(In particular, let } r_n = \frac{1}{n^2} \text{.)}
\]

Now, there exists a sequence of intervals $(a_i, b_i)$ such that:
\[
\sum_{i=1}^{k} \left| f(b_i)-f(a_i) \right| > \epsilon \quad \text{but} \quad \sum_{i=1}^{k} (b_i - a_i) < 1.
\]

take $S_1 = \bigcup_{i=1}^{k_1} (a_i, b_i)$, 

Now there exist a closed bounded interval $L_1$
such that $S_1\subset L_1$

Note that $\int_{L_1} |f'| < \infty$ since $f' \in L^1_G$

hence by continuity of measure , given 
\begin{equation}\label{r2}
\epsilon_1 =\frac{\epsilon}{2^2}\,,\,
\exists\,\delta_1 \text{ such that }\int_{E\cap L_1} |f'| < \epsilon_1 
\text{ whenever }\mu(E\cap L_1 ) <\delta_1
\end{equation}

let $m_1 = \min(\delta_1, r_1)$.

Then from (\ref{r1}), there exists $S_2$ (a finite collection of disjoint intervals)
\[
S_2 = \bigcup_{i=1}^{k_2} (a_i, b_i),
\]
such that 
\[
m(S_2) < m_1, \quad \text{and} \quad \sum_{i=1}^{k_2} \left| f(b_{i_2}) - f(a_{i_2})\right| > \epsilon.
\]
since $ f\in AC_{loc}(\mathbb{R}) \implies\int_{a_i}^{b_i} f' d\mu =  f(b_{i_2}) - f(a_{i_2})$
due to Lebesgue's fundamental theorem of calculus 

\begin{equation}\label{r3}
    \implies \sum_{i=1}^{k_2} \left| \int_{a_i}^{b_i} f' \, d\mu \right| > \epsilon.
\implies  \sum_{i=1}^{k_2}  \int_{a_i}^{b_i} \left|f' \right|\, d\mu  > \epsilon.
\implies \int_{S_2} \left|f'\right| >\epsilon
\end{equation}    
\begin{align*}
\implies&&
\int_{S_2 \setminus L_1} \left|f'\right| \, d\mu + \int_{S_2 \cap L_1} \left|f'\right| \, d\mu & > \epsilon\\
\implies && 
\int_{S_2 \setminus L_1} \left|f'\right| \, d\mu &> \epsilon - \int_{S_2 \cap L_1} \left|f'\right| \, d\mu\\ 
\implies&&& >\epsilon-\frac{\epsilon}{2^2}\quad(\because \mu(S_2\cap L_1 <\delta_1 )  \text{ from }\text{(\ref{r2}))}
\end{align*}

then consider $S_2\setminus L_1 = A_2$ and $S_1 = A_1$

In particular, after having chosen $A_n$ , consider a closed bounded interval $L_n $ such that $ \bigcup _{k = 1}^n A_k \subset L_n$, note that $\mu(L_n) <\infty $ therefore $\int_{L_n} \left|f'\right| \, d\mu < \infty$, hence by continuity of  measure , given $\epsilon_n = \frac{\epsilon}{(n+1)^2}$,there exists a $\delta_n$ such that 

\begin{equation}\label{r4}
    \int_{E\cap L_n} \left|f'\right| \, d\mu < \epsilon_n = \frac{\epsilon}{(n+1)^2} \quad \forall \mu(E\cap L_n) <\delta_n
\end{equation}
set $m_n = \min (\delta_n,r_n) $

now from \ref{r1} ,there exists a set $S_{n+1}$(which is a countable collection of disjoint intervals)
such that 
\[\mu(S_{n+1})  \leq m_n \text{  but  } \int_{S_{n+1}} \left|f'\right|\geq\epsilon
\]
\[
\implies \int_{S_{n+1}\setminus L_n} \left|f'\right|\geq\epsilon - \int_{S_{n+1}\cap L_n} \left|f'\right|\geq\epsilon - \frac{\epsilon}{(n+1)^2} \text{    from (\ref{r4})}
\]

hence choose $A_{n+1} = S_{n+1} \setminus L_n$.
Now  we have a countable disjoint collection of open sets $A_n \,,\, n\in \mathbb{N} $ such that 

\[
\mu\left(\bigcup_{n\in \mathbb{N}} A_n \right)\leq
\sum_{n\in \mathbb{N}}\mu(A_n) \leq
\sum_{n\in \mathbb{N}}r_n <\infty
\]
but 
\begin{align*}\displaystyle
    \int_{\bigcup_{n\in \mathbb{N}} A_n } \left|f'\right|d\mu&= \sum_{n = 1}^\infty \int_{ A_n } \left|f'\right|d\mu\\
    &\geq \sum_{n = 1}^\infty \epsilon- \frac{\epsilon}{n^2}\\
    &=\lim_{N\to \infty}N\epsilon - \epsilon\left(\sum_{i = 1}^N \frac{1}{i^2}\right) = \infty
\end{align*}
thereby contradicting the fact that $f' \in L^1_G(\mathbb{R})$\\

\textbf{Converse Part:}

We prove that if $f \in AC(\mathbb{R})$, then $f' \in L^1_G(\mathbb{R})$.  

$Since, f \in AC(\mathbb{R})$, hence for a given $\epsilon = 1$, $\exists \, \delta > 0$, such that if there exists a collection of disjoint intervals $(a_i, b_i)$ such that:
$\sum_{i=1}^\infty (b_i - a_i) < \delta$
then$\sum_{i=1}^\infty \left| f(a_i)-f(b_i)\right| < 1$

Now given any finite measure set $S$, there exists a set of disjoint open intervals (countable) $I_n, \, n \in \mathbb{N}$ such that
\[
S \subset \bigcup_{n=1}^\infty I_n=I, \quad \text{where} \quad \mu\left(\bigcup_{n=1}^\infty I_n\right) = \sum_{n=1}^\infty \mu(I_n) < \infty
\]
Now it suffices to show that:
\[
\int_I |f'| < \infty \text{  which will further imply}  \int_S|f'| < \infty.
\]

now since $\sum_{n=1}^\infty \mu(I_n) < \infty$,
we can have two cases:

\textbf{Case 1: $\sum_{n \in \mathbb{N}} \mu(I_n) < \delta$}

In this case choosing union over finite number of intervals $\bigcup_{k=1}^N I_k$,we can consider any arbitrary partition of this union of intervals ${P_k:k\in\{1,...,N\}}$ such that 
\[
\sum_{k=1}^N t_{a_k}^{b_k}(P_k)<1 \quad \because \sum_{k=1}^n \mu(I_k) < \delta \quad\quad I_k = (a_k,b_k)
\]
here $P_k$ represents an arbitrary partition of the $I_k$$^{th}$ interval. and $t_{a_k}^{b_k}$ represents the variation of $f$ corresponding to $P_k$. Now by taking supremum of variation with respect to all possible $P_k$'s , we have
\[
\sum_{k=1}^N T_{a_k}^{b_k}(f)<1 
\]
here $T_{a_k}^{b_k}(f) = \sup\{t_{a_k}^{b_k}(P_k) \,:\, P_k\text{ is a partition of }(a_k,b_k) \}$ 
we call $\bigcup_{k=1}^\infty I_k = I$
\[
\sum_{k=1}^N T_{a_k}^{b_k}(f)<1 \]
 , but since the function is absolutely continuous over any finite interval
\[
\therefore T_{a_k}^{b_k}(f) =  \int_{I_k} \left|f'\right| \, d\mu 
\]
\[
\implies \sum_{k=1}^N  \int_{I_k} \left|f'\right|  \, d\mu \leq 1,
\]
since this is true for all $N \in \mathbb{N}$
Therefore
\[
\sum_{k=1}^\infty \int_{I_k} |f'|d\mu \implies\int_{I} |f'| d\mu < 1.
\]

\textbf{Case 2: $\sum \mu(I_n) = M > \delta$}

let\[
\bigcup_{n=1}^\infty I_n = I.
\]
now consider $n_0 = \left\lfloor \frac{M}{\frac{\delta}{2}} \right\rfloor$, now note that $M - n_0 \frac{\delta}{2} < \frac{\delta}{2}$.  
also there exist a large finite $N_0 \in \mathbb{N}$, such that:
\[
n_0\frac{\delta}{2} \leq\sum_{n=1}^{N_0} m(I_n) \leq M .
\]

now each of these $N_0$ number of intervals can be subdivided into intervals or union of intervals having length in between $\delta$ and $\frac{\delta}{2},$  let $N'$ be the number of such new disjoint collection of intervals $J_n$, with 
$\frac{\delta}{2}<\mu(J_n)<\delta$ for each $n$.

hence
\[
\int_{J_n} |f'| = \text{Total Variation}(J_n) \leq 1, \quad \forall n\in \mathbb{N}.
\]

Now
\begin{equation}\label{n1}
\int_{\bigcup_{n=1}^{N'} J_n} |f'| d\mu\leq N'\leq n_0
\end{equation}
also $I-\bigcup_{n=1}^{N'} J_n$ has measure $<$ $\frac{\delta}{2}$
hence by proof of case-I
\begin{equation}\label{n2}
\int_{I\setminus\bigcup_{n=1}^{N'} J_n} |f'| d\mu\leq 1
\end{equation}

therefore from (\ref{n1}) and (\ref{n2}) , we get 
\[
\therefore \int_{I} |f'| d\mu\leq n_0+1
\]
\[
\therefore \int_{S} |f'| d\mu\leq n_0+1
\]
hence finite.
\\
Without loss of generality, we can conclude that if  $f\in AC_{loc}(\mathbb{I})$ and $f' \in L^1_G(\mathbb{I}) \iff f\in AC(\mathbb{I}) $ where $\mathbb{I}$ is an interval in $\mathbb{R}$ .
\end{proof}

\begin{remark}
Consider the set 
\[
A = \bigcup_{n=1}^\infty\left[2n,2n+\frac{1}{n^2}\right]
\]
Now consider the function $f$ in example \ref{ex1} which has been introduced earlier as an example of a function that is locally absolutely continuous but not globally absolutely continuous.
Note that A is a finite measure set but $\int_A|f'|d\mu = \infty$ hence $f' \notin L^1_G$ therefore $f\notin AC(\mathbb{R})$.
Thus, this theorem provides us with a simpler criterion for checking whether a function is absolutely continuous over $\mathbb{R}$ or not
\end{remark}

\begin{theorem}
$L^1_G(\mathbb{I})$ is a proper subset of $ L^1_H(\mathbb{I}) $ where,

$L^1_G(\mathbb{R})$ = $\{$ space of measurable functions $ f:\mathbb{R} \to \mathbb{R} $ such that $ \int_K |f| < \infty$ where K is a lebesgue measurable subset of I and with finite measure.$\}$ and 

$L^1_H(\mathbb{R})$ = $\{$set of all measurable functions $ f : \mathbb{R} \to \mathbb{R} |$  for a given $ f $ there exists $ M > 0 : \mu\{x \in \mathbb{I} \mid |f(x)| \geq M\} <\infty\} $.
\end{theorem}

\begin{proof}
We show that $\left(L^1_H(\mathbb{I})\right)^c \subset \left(L^1_G(\mathbb{I})\right)^c$ \\

Let $f \in \left(L^1_H(\mathbb{I})\right)^c$. \\

Given $n \in \mathbb{N}$, define 
\[
S_n = \{x \in \mathbb{I} \mid |f(x)| \geq n\}.
\]
Now $\mu(S_n) = \infty \quad \forall n \in \mathbb{N}$.

Define:
\[
F_{(m)_n} = \left[ -\frac{m}{2n^2}, \frac{m}{2n^2} \right], \quad m, n \in \mathbb{N}.
\]

Now, there exists a sequence of sets:
\[
F_{(a_1)_1} \subset F_{(a_2)_2} \subset F_{(a_3)_3} \subset \cdots \subset F_{(a_n)_n}\subset\cdots
\]
where $(a_{n})<(a_{n+1})$.
such that
\[
 1 \leq \mu\left(F_{(a_1)_1} \cap S_1\right) \leq 2.
\]
\[
\frac{1}{2^2} \leq \mu\left(F_{(a_2)_2}\setminus F_{(a_1)_1} \cap S_2\right)\leq \frac{2}{2^2}
\]
In general,
\[
\frac{1}{n^2} \leq \mu\left(F_{(a_n)_n}\setminus F_{(a_{n-1})_{n-1}} \cap S_n\right)\leq \frac{2}{n^2}
\]

Define $F_{(a_0)_0} = \phi$.

Now, consider the set:
\[
G = \bigcup_{n=1}^\infty G_n,
\]
where
\[
G_n = \left(F_{(a_n)_n} \setminus F_{(a_{n-1})_{n-1}} \cap S_n\right).
\]
Note that the sets $G_n$ are disjoint, \[\therefore\mu(G_n) = \sum_{n=1}^\infty \mu(G_n) < \sum_{n=1}^\infty \frac{2}{n^2} < \infty\]

However,
\[
\int_G |f| \, d\mu \geq \sum_{n=1}^\infty n \times \frac{1}{n^2} = \sum_{n=1}^\infty \frac{1}{n} = \infty.
\]

\end{proof}

\begin{remark}
Consider the function $f$ mentioned in example \ref{ex1}. Given any real number $r>0$, there exist an interval of some finite length around each of the even integral points, such that $|f'(x)|>r$ in that interval hence for any $r>0, \mu\{x\in \mathbb{R}\mid |f'(x)|>r \} = \infty$ therefore 
\[
f'\notin L^1_H(\mathbb{R})\implies f'\notin L^1_G(\mathbb{R})\implies f\notin AC(\mathbb{R})
\]
\end{remark}

\end{document}